\documentclass{article}
\usepackage{spconf}

\def\bsg{{\boldsymbol{g}}}

\def\bsx{{\boldsymbol{x}}}

\def\bsH{{\boldsymbol{H}}}

\def\bsK{{\boldsymbol{K}}}
\def\bsL{{\boldsymbol{L}}}

\def\bsM{{\boldsymbol{M}}}

\def\bsQ{{\boldsymbol{Q}}}

\usepackage{verbatim}
\usepackage{float}
\usepackage{bm}
\usepackage[cmex10]{amsmath}
\usepackage{amssymb}

\usepackage{graphicx}
\graphicspath{{./figures/}}

\def\sgn{\mathop{\rm sgn}}
\usepackage{amsthm}

\usepackage{graphicx}
\usepackage{subcaption}
\usepackage{url}
\usepackage{breakurl}
\usepackage{tikz}
\usepackage{verbatim}
\usepackage{forest}
\usepackage[T1]{fontenc}
\usepackage{pifont}
\usepackage{mathtools}
\usetikzlibrary{shapes,arrows,positioning,calc}
\usepackage{cite}
\usepackage[english]{babel}
\usepackage[utf8]{inputenc}

\usepackage{algorithm,algorithmic}

\newtheorem{theorem}{Theorem}
\newtheorem{assumption}{Assumption}

\newtheorem{lemma}{Lemma}

\newtheorem{remark}{Remark}

\usepackage{color}

\makeatletter

\tikzstyle{block} = [draw, thick, rectangle, 
minimum height=3em, minimum width=6em]
\tikzstyle{pt} = [coordinate]
\tikzset{
	block/.style={
		draw, 
		rectangle, 
		minimum height=0.8cm, 
		minimum width=0.8cm, 
		align=center
	}
}

\newcommand{\diag}{{\rm diag}}

\DeclareMathOperator{\Deg}{Deg}

\DeclareMathOperator{\col}{col}

\usetikzlibrary{arrows,shapes,positioning,shadows,trees}

\tikzset{
  basic/.style  = {draw, text width=5cm, drop shadow, font=\sffamily, rectangle},
  root/.style   = {basic, rounded corners=2pt, thin, align=center,
                   fill=yellow!60},
  level 2/.style = {basic, rounded corners=6pt, thin,align=center, fill=yellow!30,
                   text width=8em},
  level 3/.style = {basic, thin, align=left, fill=pink!40, text width=3cm}
}

\usepackage{adjustbox}
\usepackage{booktabs} 
\usepackage{makecell}
\allowdisplaybreaks
\begin{document}

\title{Zeroth-Order Stochastic Coordinate Methods for Decentralized Non-convex Optimization}

\name{Shengjun~Zhang$^\dagger$, Tan~Shen$^\ddagger$, Hongwei~Sun$^\ddagger$, Yunlong~Dong$^\mathsection$, Dong~Xie$^{**}$, and Heng~Zhang$^*$}
\address{$^\dagger$AlgoRhythm Inc.\\$^\ddagger$School of Artificial Intelligence and Automation, Huazhong University of Science and Technology\\
$^\mathsection$Department of Automation, Tsinghua University\\
$^{**}$Lenovo Research
\\$^*$School of Electronic Information and Electrical Engineering, Shanghai Jiao Tong University}
%
%
%

\maketitle

\begin{abstract}            

In this paper, we first propose a \underline{Z}eroth-\underline{O}rder c\underline{O}ordinate \underline{M}ethod~(ZOOM) to solve the stochastic optimization problem over a decentralized network with only zeroth-order~(ZO) oracle feedback available.
Moreover, we equip a simple mechanism "powerball" to ZOOM and propose ZOOM-PB to accelerate the convergence of ZOOM.
Compared with the existing methods, we verify the proposed algorithms through two benchmark examples in the literature, namely the black-box binary classification and the generating adversarial examples from black-box DNNs in order to compare with the existing state-of-the-art centralized and distributed ZO algorithms.
The numerical results demonstrate a faster convergence rate of the proposed algorithms.


\end{abstract}


\section{Introduction\label{sec:Introduction}}

Decentralized optimization problems have been investigated for a long time and can be traced back to the early 1980's~\cite{tsitsiklis1984problems}, where the parallel and distributed computation was first introduced.
Over the last decade, along with the rapid developments of communication, signal processing and control for networked multi-agent systems, decentralized optimization has been studied and applied in various many areas for instance, power systems~\cite{luo2005review}, machine learning~\cite{choi2009distributed}, swarm robotics control~\cite{jaleel2020distributed}, for the recent review and progress, please refer to the survey~\cite{yang2019survey}.

Generally speaking, in a decentralized optimization problem, a certain number of agents~(nodes) cooperatively find the optimal solution to Eq.~\eqref{zerosg:eqn:xopt:org} by exchanging their own local information only under a network topology, which is formulated as

\begin{align}\label{zerosg:eqn:xopt:org}
 \min_{x\in \mathbb{R}^p} f(x)=\frac{1}{n}\sum_{i=1}^n f_{i}(x),
\end{align}
where $n$ is the number of agents, $x\in \mathbb{R}^p$ denotes the optimization variable with dimension $p$, $f_i: \mathbb{R}^{p}\rightarrow \mathbb{R}$ is the local smooth (possibly non-convex) cost function of agent $i$.

To address the problem in Eq.~\eqref{zerosg:eqn:xopt:org}, many algorithms have been proposed based on the gradient~(first-order) information and ADMM, e.g.~\cite{nedic2009distributed, shi2015extra, chang2016proximal}.
Moreover, with the growth of the number of local samples, computing the full gradient of the local cost function becomes more difficult.
As a result, algorithms based on Stochastic Gradient Descent~(SGD) are more favorable in decentralized machine learning and deep learning applications~\cite{wang2021distributed}.

Unfortunately, in many scenarios, the deceptively simple gradient information is not available or too expensive to compute \cite{conn2009introduction,audet2017derivative,larson2019derivative}.
For instance, the simulation based optimization problems~\cite{spall2005introduction},  the universal attacking problems of deep neural networks~\cite{goodfellow2014explaining, liu2018zeroth, chen2019zo}, data generating processes problems~\cite{chen2017zoo}. 
In recent years, distributed zeroth-order~(ZO) optimization problems have gained more and more attention and have been applied into networked agent systems, e.g~\cite{yuan2014randomized,sahu2018distributed,
pang2019randomized,yu2019distributed,tang2020distributedzero, beznosikov2019derivative, yi2019linearTAC,hajinezhad2019zone,yi2021zerothorder,zhang2021convergence}.

Most of the ZO algorithms aforementioned consider the deterministic decentralized optimization problems, which are formulated exactly as Eq.~\eqref{zerosg:eqn:xopt:org}.
On the contrary, some applications require that the data samples are received by the local agents in a streaming way, for example the online learning problems.
Moreover, in the era of big data, with the more data people can obtain and the higher dimension of data people can measure, a deterministic algorithm will lead to higher computational burdens and demands, which motivates us to focus on stochastic problem settings. 
Such problems can be formulated into the stochastic distributed optimization problem
\begin{align}\label{zo_pb:eqn:xopt}
 \min_{x\in \mathbb{R}^p} f(x)=\frac{1}{n}\sum_{i=1}^n\mathbb{E}_{\xi_i}[F_i(x,\xi_i)],
\end{align}
where $\xi_i$ is a random variable with dimension $p$, and $F_i(\cdot,\xi_i): \mathbb{R}^{p}\rightarrow \mathbb{R}$ is the stochastic function.

Under the stochastic distributed settings in the exact form of~\eqref{zo_pb:eqn:xopt}, only a few works \cite{hajinezhad2019zone,yi2021zerothorder,zhang2021convergence} exist in the literature. 
ZONE-M in \cite{hajinezhad2019zone} achieves the convergence rate of $\mathcal{O}(p^2n/T)$ with a very high sampling size of $\mathcal{O}(T)$ per iteration.
ZODPDA~\cite{yi2021zerothorder} and ZODIAC~\cite{zhang2021convergence} achieve the convergence rate of $\mathcal{O}(\sqrt{p}/\sqrt{nT})$ with $2n$ points sampled per iteration, which are more suitable for high dimensional decision variables cases.

In this paper, to address the stochastic ZO distributed problems in the form of Eq.~\eqref{zo_pb:eqn:xopt}, we first propose a \underline{Z}eroth-\underline{O}rder c\underline{O}ordinate \underline{M}ethod~(ZOOM) to solve them.
To accelerate the convergence of ZOOM, we apply a simple but efficient method namely "powerball" to ZOOM and propose ZOOM-PB, an acceleration variation of ZOOM.
Unlike \cite{hajinezhad2019zone, yi2021zerothorder, zhang2021convergence}, the proposed algorithms are simpler to implement but have similar convergence results.
Extensive numerical examples are provided to demonstrate the efficacy of the proposed algorithms and compare with existing algorithms through benchmark examples in the ZO optimization literature.

\section{Background and Assumptions}\label{zo_pb:sec-preliminary}
The following section discusses background of graph theory, smooth functions, the gradient estimator we choose, and additional assumptions used in this paper.

\subsection{Graph Theory}
Agents communicate with their neighbors through an underlying network, which is modeled by an undirected graph $\mathcal G=(\mathcal V,\mathcal E)$, where $\mathcal V =\{1,\dots,n\}$ is the agent set, $\mathcal E
\subseteq \mathcal V \times \mathcal V$ is the edge set, and $(i,j)\in \mathcal E$ if agents $i$ and $j$ can communicate with each other.
For an undirected graph $\mathcal G=(\mathcal V,\mathcal E)$, let $\mathcal{A}=(a_{ij})$ be the associated weighted adjacency matrix with $a_{ij}>0$ if $(i,j)\in \mathcal E$ and zero otherwise.
It is assumed that $a_{ii}=0$ for all $i\in [n]$. 
Let $\deg_i=\sum\limits_{j=1}^{n}a_{ij}$ denotes the weighted degree of vertex $i$, then the degree matrix of graph $\mathcal G$ is $\Deg=\diag([\deg_1, \cdots, \deg_n])$. and the Laplacian matrix associated with $\mathcal G$ is $L=(L_{ij})=\Deg-\mathcal{A}$. 
Moreover, we use $\rho(\cdot)$ to describe the spectral radius for matrices and $\rho_2(\cdot)$ to indicate the minimum
positive eigenvalue for matrices that have positive eigenvalues.
Additionally, we denote $K_n={\bm I}_n-\frac{1}{n}{\bm 1}_n{\bm 1}^{\top}_n$, $\bsL=L\otimes {\bm I}_p$, $\bsK=K_n\otimes {\bm I}_p$, $\bsH=\frac{1}{n}({\bm 1}_n{\bm 1}_n^\top\otimes{\bm I}_p)$. Moreover, from Lemmas~1 and 2 in \cite{Yi2018distributed}, we know there exists an orthogonal matrix $[r \ R]\in \mathbb{R}^{n \times n}$ with $r=\frac{1}{\sqrt{n}}\mathbf{1}_n$ and $R \in \mathbb{R}^{n\times (n-1)}$ such that $R\Lambda_1^{-1}R^{\top}L=LR\Lambda_1^{-1}R^{\top}=K_n$, and $\frac{1}{\rho(L)}K_n\leq R\Lambda_1^{-1}R^{\top}\le\frac{1}{\rho_2(L)}K_n$, where $\Lambda_1=\diag([\lambda_2,\dots,\lambda_n])$ with $0<\lambda_2\leq\dots\leq\lambda_n$ being the eigenvalues of the Laplacian matrix $L$.

\subsection{Smooth Function}
A function $f(x):~\mathbb{R}^p\mapsto\mathbb{R}$ is smooth with constant $L_f>0$ if it is differentiable and
\begin{align}\label{zo_pb:nonconvex:smooth}
\|\nabla f(x)-\nabla f(y)\|\le L_{f}\|x-y\|,~\forall x,y\in \mathbb{R}^p.
\end{align}

\subsection{Gradient Approximation}

Denote a random subset of the coordinates $\mathcal{S} \subseteq \{1, 2, \dots, p\}$ where the cardinality of $\mathcal{S}$ is $|\mathcal{S}| = n_{c}$. 
We provide two options of gradient approximation, denoted $g^e_{i}$ and defined by~\eqref{zo_pb:gradient:model2-st} and~\eqref{zo_pb:gradient:model2-st2}.

\begin{align}
&g^e_{i}=\frac{p}{n_{c}}\sum_{i\in \mathcal{S}}\frac{(F(x+\delta_{i} e_{i},\xi)-F(x,\xi))}{\delta_{i}}e_{i}
\label{zo_pb:gradient:model2-st}
\end{align}

\begin{align}
&g^e_{i}=\frac{p}{n_{c}}\sum_{i\in \mathcal{S}}\frac{(F(x+\delta_{i} e_{i},\xi)-F(x-\delta_{i} e_{i},\xi))}{2\delta_{i}}e_{i}
\label{zo_pb:gradient:model2-st2}
\end{align}


\subsection{Powerball Function}

Define the function 
\begin{equation}\label{zo_pb:pb:def}
\sigma(x, \gamma) = \sgn(x) |x|^{\gamma}
\end{equation}
where $\gamma \in [\frac{1}{2}, 1]$.
Note that when $\gamma = 1$, $\sigma(x, 1)$ reduces to $x$.
\subsection{Assumptions}

\begin{assumption}\label{zo_pb:ass:graph}
The undirected graph $\mathcal G$ is connected.
\end{assumption}

\begin{assumption}\label{zo_pb:ass:optset}
The optimal set $\mathbb{X}^*$ is nonempty and the optimal value $f^*>-\infty$.
\end{assumption}

\begin{assumption}\label{zo_pb:ass:zeroth-smooth}
For almost all $\xi_i$, the stochastic ZO oracle $F_i(\cdot,\xi_i)$ is smooth with constant $L_f>0$.
\end{assumption}

\begin{assumption}\label{zo_pb:ass:zeroth-variance}
The stochastic gradient $\nabla_xF_i(x,\xi_i)$ has bounded variance for any $j$th coordinate of $x$, i.e., there exists $\zeta\in\mathbb{R}$ such that $\mathbb{E}_{\xi_i}[(\nabla_xF_i(x,\xi_i)-\nabla f_i(x))_{j}^2]\le\zeta^2,~\forall i\in[n],~\forall j\in[p],~\forall x\in\mathbb{R}^p$. It also implies that $\mathbb{E}_{\xi_i}[\|\nabla_xF_i(x,\xi_i)-\nabla f_i(x)\|^2]\le\sigma^2_1\triangleq p \zeta^2,~\forall i\in[n],~\forall x\in\mathbb{R}^p$.
\end{assumption}

\begin{assumption}\label{zo_pb:ass:fig}
Local cost functions are similar, i.e.,
there exists $\sigma_2\in\mathbb{R}$ such that $\|\nabla f_i(x)-\nabla f(x)\|^2\le\sigma^2_2,~\forall i\in[n],~\forall x\in\mathbb{R}^p$.
\end{assumption}


\section{Proposed Algorithms}

\subsection{Algorithm Description}
We summarize the proposed algorithms ZOOM and ZOOM-PB as Algorithm~\ref{zerosc-p:algorithm-random} and Algorithm~\ref{zerosc-p:algorithm-random-pb} respectively.

\begin{algorithm}[ht!]
\caption{ZOOM}
\label{zerosc-p:algorithm-random}
\begin{algorithmic}[1]
\STATE \textbf{Input}: positive constant $\alpha$, $\eta$ and $\{\delta_{i,k}\}$.
\STATE \textbf{Initialize}: $ x_{i,0}\in\mathbb{R}^p,~\forall i\in[n]$.
\FOR{$k=0,1,\dots$}
\FOR{$i=1,\dots,n$  in parallel}
\STATE  Broadcast $x_{i,k}$ to $\mathcal{N}_i$ and receive $x_{j,k}$ from $j\in\mathcal{N}_i$;
\STATE \textbf{Option 1:} \\sample $F_i(x_{i,k}+\delta_{i,k}e_{i,k},\xi_{i,k})$, $F_i(x_{i,k},\xi_{i,k})$, and update $g^e_{i,k}$ by \eqref{zo_pb:gradient:model2-st};
\STATE \textbf{Option 2:} \\sample  $F_i(x_{i,k}+\delta_{i,k}e_{i,k},\xi_{i,k})$, $F_i(x_{i,k}-\delta_{i,k}e_{i,k},\xi_{i,k})$, and update $g^e_{i,k}$ by \eqref{zo_pb:gradient:model2-st2};
\STATE  Update $x_{i,k+1} =x_{i,k}-\alpha\sum_{j=1}^nL_{ij}x_{j,k}-\eta g^e_{i,k}$.
\ENDFOR
\ENDFOR
\STATE  \textbf{Output}: $\{\mathbf{x_{k}}\}$.
\end{algorithmic}
\end{algorithm}

\begin{algorithm}[ht!]
\caption{ZOOM-PB}
\label{zerosc-p:algorithm-random-pb}
\begin{algorithmic}[1]
\STATE \textbf{Input}: positive constants $\alpha$, $\gamma$, $\eta$ and $\{\delta_{i,k}\}$.
\STATE \textbf{Initialize}: $ x_{i,0}\in\mathbb{R}^p,~\forall i\in[n]$.
\FOR{$k=0,1,\dots$}
\FOR{$i=1,\dots,n$  in parallel}
\STATE  Broadcast $x_{i,k}$ to $\mathcal{N}_i$ and receive $x_{j,k}$ from $j\in\mathcal{N}_i$;
\STATE \textbf{Option 1:} \\sample $F_i(x_{i,k}+\delta_{i,k}e_{i,k},\xi_{i,k})$ $F_i(x_{i,k},\xi_{i,k})$, and update $g^e_{i,k}$ by \eqref{zo_pb:gradient:model2-st};
\STATE \textbf{Option 2:} \\sample  $F_i(x_{i,k}+\delta_{i,k}e_{i,k},\xi_{i,k})$, $F_i(x_{i,k}-\delta_{i,k}e_{i,k},\xi_{i,k})$, and update $g^e_{i,k}$ by \eqref{zo_pb:gradient:model2-st2};
\STATE  Update \\ $x_{i,k+1} =x_{i,k}-\alpha\sum_{j=1}^nL_{ij}x_{j,k}-\eta \sigma(g^e_{i,k}, \gamma)$.\label{zo_pb:alg:update}
\ENDFOR
\ENDFOR
\STATE  \textbf{Output}: $\{\mathbf{x_{k}}\}$.
\end{algorithmic}
\end{algorithm}

\subsection{Convergence Analysis}

\begin{theorem}\label{zerosg-p:thm-sg-smT-ZOOM}
Suppose Assumptions~\ref{zo_pb:ass:graph}--\ref{zo_pb:ass:fig} hold. For any given $T\ge n^3/p$, let $\{\bsx_k,k=0,\dots,T\}$ be the output generated by Algorithm~\ref{zerosc-p:algorithm-random} with
\begin{align}\label{zerosg-p:step:eta2-sm}
&\alpha\in\Big(0,\frac{\rho_2(L)}{2\rho(L^2)}\Big),~\eta=\frac{\sqrt{n}}{\sqrt{pT}},\nonumber\\
&\delta_{i,k}\le\frac{\kappa_\delta}{p^{\frac{1}{4}}n^{\frac{1}{4}}(k+1)^{\frac{1}{4}}},~\forall k\le T,
\end{align}
where $\kappa_\delta>0$, then
\begin{subequations}
\begin{align}
&\frac{1}{T}\sum_{k=0}^{T-1}\mathbf{E}[\|\nabla f(\bar{x}_k)\|^2]
=\mathcal{O}(\frac{\sqrt{p}}{\sqrt{nT}})+\mathcal{O}(\frac{n}{T}),\label{zerosg-p:thm-sg-sm-equ3-ZOOM}\\
&\mathbf{E}[f(\bar{x}_{T})]-f^*=\mathcal{O}(1),\label{zerosg-p:thm-sg-sm-equ4-ZOOM}\\
&\frac{1}{T}\sum_{k=0}^{T-1}\mathbf{E}\Big[\frac{1}{n}\sum_{i=1}^{n}
\|x_{i,k}-\bar{x}_k\|^2\Big]=\mathcal{O}(\frac{n}{T}).\label{zerosg-p:thm-sg-sm-equ3.1-ZOOM}
\end{align}
\end{subequations}
\end{theorem}

\begin{proof}
Please see Remark~\ref{rmk:th1th2}.
\end{proof}


\begin{theorem}\label{zerosg-p:thm-sg-smT-ZOOM-PB}
Suppose  Assumptions~\ref{zo_pb:ass:graph}--\ref{zo_pb:ass:fig} hold. For any given $T\ge n^3/p$, let $\{\bsx_k,k=0,\dots,T\}$ be the output generated by Algorithm~\ref{zerosc-p:algorithm-random-pb} with
\begin{align}\label{zerosg-p:step:eta2-sm-ZOOM-PB}
&\alpha\in\Big(0,\frac{\rho_2(L)}{2\rho(L^2)}\Big),~\eta=\frac{\sqrt{n}}{\sqrt{pT}},~\gamma \in [\frac{1}{2}, 1],\nonumber\\
&\delta_{i,k}\le\frac{\kappa_\delta}{p^{\frac{1}{4}}n^{\frac{1}{4}}(k+1)^{\frac{1}{4}}},~\forall k\le T, 
\end{align}
where $\kappa_\delta>0$, then

\begin{subequations}
\begin{align}
&\frac{1}{T}\sum_{k=0}^{T-1}\mathbf{E}[\|\nabla f(\bar{x}_k)\|_{1+\gamma}^2]
=\mathcal{O}(\frac{\sqrt{p}}{\sqrt{nT}})+\mathcal{O}(\frac{n}{T}),\label{zerosg-p:thm-sg-sm-equ3-ZOOM-PB}\\
&\mathbf{E}[f(\bar{x}_{T})]-f^*=\mathcal{O}(1),\label{zerosg-p:thm-sg-sm-equ4-ZOOM-PB}\\
&\frac{1}{T}\sum_{k=0}^{T-1}\mathbf{E}\Big[\frac{1}{n}\sum_{i=1}^{n}
\|x_{i,k}-\bar{x}_k\|^2\Big]=\mathcal{O}(\frac{n}{T}).\label{zerosg-p:thm-sg-sm-equ3.1-ZOOM-PB}
\end{align}
\end{subequations}
\end{theorem}

\begin{remark}\label{rmk:th1th2}
We can see that ZOOM is a special case of ZOOM-PB when $\gamma = 1$, therefore, we omit the proof of Theorem~\ref{zerosg-p:thm-sg-smT-ZOOM}.
\end{remark}

In order to prove Theorem~\ref{zerosc-p:algorithm-random-pb}, we introduce the following lemmas.

\begin{lemma}(Lemma 2 in~\cite{zhang2021convergence})\label{zo_pb:lemma:variance}
Consider $f(x) = \mathbb{E}_{\xi} [F(x, \xi)]$, we have the following relationship,
\begin{align}
&\mathbb{E} \Big[ \| g^e_{i} \|  ^2\Big] \nonumber \\
&\leq 2(p-1)\left\| \nabla f (x) \right\|^2 + 2p \sigma^2_1 + \frac{3p^{2}}{n_{c}} \left( \zeta^2 + \frac{L_{f}^2 \delta_{k}^2}{2} \right) \nonumber \\
&\quad+ \frac{p^{2}L_{f}^2 \delta_{k}^2}{2} \label{zo_pb:eq_bd_var_CGE}
\end{align}
\end{lemma}
where $\delta_{k} = \max\{ \delta_{i} \}, i \in [p]$.

\begin{lemma}(Lemma 2 in~\cite{zhang2021accelerated})\label{zo_pb:lemma:pb}
By using the powerball term in~\eqref{zo_pb:pb:def} and when $\gamma \in [\frac{1}{2}, 1]$, we have $\Big\Vert\sigma(g^e_{i}, \gamma)\Big\Vert^2 \leq \Big\Vert g^e_{i} \Big\Vert^{2}_{1+\gamma}$.
\end{lemma}

\begin{lemma}(Lemma 3 in~\cite{zhang2021accelerated})\label{zo_pb:zerosg:lemma:grad-st}
Suppose Assumptions~~\ref{zo_pb:ass:zeroth-smooth}--~\ref{zo_pb:ass:fig} hold. Let $\{\bsx_k\}$ be the sequence generated by Algorithm~\ref{zerosc-p:algorithm-random-pb}, $\bsg^e_k=\col(g^e_{1,k},\dots,g^e_{n,k})$, $\bsg^0_k=n\nabla{f}(\bar{\bsx}_k)$, $\bar{\bsg}_k^0=\bsH\bsg^0_{k}={\bm 1}_n\otimes\nabla f(\bar{x}_k)$, then
\begin{subequations}
\begin{align}
\mathbb{E}\Big[\|\bsg^e_k\|_{1+\gamma}^2\Big]
&\le  6(p-1)\|\bar{\bsg}_{k}^0\|_{1+\gamma}^2+6(p-1)L_f^2\|\bsx_{k}\|^2_{\bsK}\nonumber\\
&\quad+6n(p-1)\sigma^2_2 + \frac{3np^{2}}{n_{c}} \left( \zeta^2 + \frac{L_{f}^2 \delta_{k}^2}{2} \right)\nonumber \\
&\quad+2np \sigma^2_1+ \frac{np^{2}L_{f}^2 \delta_{k}^2}{2}\label{zo_pb:zerosg:rand-grad-esti2}\\
\|\bsg^0_{k+1}\|^2&\le 3(\eta^2L_f^2\|\bsg^e_{k}\|^2+n\sigma^2_2
+\|\bar{\bsg}_{k}^0\|^2).\label{zo_pb:zerosg:rand-grad-esti4}
\end{align}
\end{subequations}
\end{lemma}
%

\begin{lemma}\label{zo_pb:lemma:sg2-T}
Suppose Assumptions~\ref{zo_pb:ass:graph}--\ref{zo_pb:ass:fig} hold, and we have fixed parameters  $\alpha \in (0, \frac{\rho_2(L)}{2\rho(L^2)})$, and \\
$\eta \in (0,\min \{\frac{2\alpha\rho_2(L)-4\alpha^2\rho(L^2)}{9L_f^{2}}, \frac{\alpha\rho_2(L)}{48p[(1+2\alpha\rho_2(L))+\alpha\rho_2(L)L_f]}\}\big]$
are constants. Let $\{\bsx_k\}$ be the sequence generated by Algorithm~\ref{zerosc-p:algorithm-random-pb}, then
\begin{subequations}
\begin{align}
\mathbb{E}[W_{k+1}]
&\le   W_{k}-\frac{\epsilon}{2}\|\bsx_k\|^2_{\bsK}
-\frac{1}{8}\eta\|\bar{\bsg}^0_{k}\|_{1+\gamma}^2\nonumber\\
&\quad+\mathcal{O}(np)\eta^2 + \mathcal{O}(np^2)\eta\delta_k^2,
\label{zo_pb:zerosg:sgproof-vkLya2T}\\
\mathbb{E}[W_{4,k+1}]
&\le  W_{4,k}+2\eta L_f^2\|\bsx_k\|^2_{\bsK}-\frac{1}{8}\eta\|\bar{\bsg}_{k}^0\|^2\nonumber\\
&\quad+\mathcal{O}(p)\eta^2
+\mathcal{O}(np)\eta\delta^2_k.\label{zo_pb:zerosg:v4kspeed}
\end{align}
\end{subequations}
\end{lemma}

\begin{proof}
We provide the proof of Lemma~\ref{zo_pb:lemma:sg2-T} in the appendix.
\end{proof}

We are now ready to prove Theorem~\ref{zerosg-p:thm-sg-smT-ZOOM-PB}.

Denote
\begin{align*}
\hat{V}_k=\|\bm{x}_k\|^2_{\bsK}+n(f(\bar{x}_k)-f^*).
\end{align*}
We have
\begin{align}
&W_{k}\nonumber\\
&=\frac{1}{2}\|\bsx_{k}\|^2_{\bsK}+n(f(\bar{x}_k)-f^*)\nonumber\\
&\ge\frac{1}{2}\|\bsx_{k}\|^2_{\bsK}
+\frac{1}{2}\Big(\frac{1}{\rho(L)}+\kappa_1\Big)
\Big\|{\beta_k}\bsg_k^0\Big\|^2_{\bsK}\nonumber\\
&~~~-\frac{1}{2\kappa_1}\|\bsx_{k}\|^2_{\bsK}
-\frac{1}{2}\kappa_1\Big\|\bsg_k^0\Big\|^2_{\bsK}
+n(f(\bar{x}_k)-f^*)\nonumber\\
&\ge\min\Big\{\frac{1}{2\rho(L)},~\frac{\kappa_1-1}{2\kappa_1}\Big\}\hat{V}_k\ge0,\label{zo_pb:zerosg:vkLya3}
\end{align}
where $\kappa_1$ is given in Appendix~\ref{zo_pb:proof:lemma4}.

From $\eta=\sqrt{n}/\sqrt{pT}$ and $T\ge n^3/p$, we know that all conditions needed in Lemma~\ref{zo_pb:lemma:sg2-T} are satisfied. So \eqref{zo_pb:zerosg:sgproof-vkLya2T} and \eqref{zo_pb:zerosg:v4kspeed} hold.

From \eqref{zo_pb:zerosg:sgproof-vkLya2T}, $\eta=\sqrt{n}/\sqrt{pT}$, and $\delta_{i,k}\le\kappa_\delta/(pn(k+1))^{1/4}$ as stated in \eqref{zerosg-p:step:eta2-sm}. Summing~\eqref{zo_pb:zerosg:sgproof-vkLya2T} over $k \in [0, T]$ and applying~\eqref{zo_pb:zerosg:vkLya3}, we have
\begin{align}
&\frac{1}{T+1}\sum_{k=0}^{T}\mathbb{E}[\frac{1}{n}\sum_{i=1}^{n}\|x_{i,k}-\bar{x}_k\|^2]\nonumber\\
&\le\frac{1}{\kappa_4}\Big(\frac{W_{0}}{n(T+1)}
+\frac{\mathcal{O}(n)\eta\delta_k^2}{T}
+\frac{\mathcal{O}(n/p)\eta^2\kappa_\delta}{\sqrt{T(T+1)}}\Big) \nonumber\\
&=\mathcal{O}(\frac{n}{T}),
\label{zo_pb:zerosg:thm-sg-sm-equ3.1p}
\end{align}
where $W_{0} = \mathcal{O}(n)$, $\frac{W_{0}}{n(T+1)}=\mathcal{O}(\frac{1}{T})$, $\frac{n\mathcal{O}(p^2)\eta\delta_k^2}{T} = \mathcal{O}(\frac{n}{T})$, and $\frac{\mathcal{O}(n/p)\eta^2\kappa_\delta}{\sqrt{T(T+1)}}=\mathcal{O}(\frac{n}{pT})$ , which gives~\eqref{zerosg-p:thm-sg-sm-equ3.1-ZOOM-PB}.

From~\eqref{zo_pb:zerosg:v4kspeed},~\eqref{zerosg-p:thm-sg-sm-equ3.1-ZOOM-PB}, and~\eqref{zo_pb:zerosg:vkLya3}, summing~\eqref{zo_pb:zerosg:v4kspeed} over $k \in [0, T]$ similar to the way to get~\eqref{zerosg-p:thm-sg-sm-equ3.1-ZOOM-PB}, we have
\begin{align}\label{zo_pb:zerosg:thm-sg-sm-equ3p}
&\frac{1}{T+1}\sum_{k=0}^{T}\mathbb{E}[\|\nabla f(\bar{x}_k)\|_{1+\gamma}^2]=\frac{1}{n(T+1)}\sum_{k=0}^{T}\mathbb{E}[\|\bar{\bsg}_{k}^0\|_{1+\gamma}^2]\nonumber\\
&\le 8\Big(\frac{W_{4,0}}{n(T+1)\eta}
+\frac{2L_f^2}{n(T+1)}\sum_{k=0}^{T}\mathbb{E}[\|\bsx_k\|^2_{\bsK}]+\frac{\mathcal{O}(p)}{n}\nonumber\\
&~~~+\frac{\mathcal{O}(\sqrt{np})}{n\sqrt{T+1}}\Big).
\end{align}
Noting that $\eta=\kappa_2/\beta_k=\sqrt{n}/\sqrt{pT}$, and $n/T<\sqrt{p}/\sqrt{nT}$ due to $T> n^3/p$, from~\eqref{zo_pb:zerosg:thm-sg-sm-equ3p} and~\eqref{zo_pb:zerosg:thm-sg-sm-equ3.1p}, we have
\begin{align*}
\frac{1}{T}\sum_{k=0}^{T-1}\mathbb{E}[\|\nabla f(\bar{x}_k)\|_{1+\gamma}^2]
&=\mathcal{O}(\frac{\sqrt{p}}{\sqrt{nT}})
+\mathcal{O}(\frac{n}{T}),
\end{align*}
which gives~\eqref{zerosg-p:thm-sg-sm-equ3-ZOOM-PB}.

Summing~\eqref{zo_pb:zerosg:v4kspeed} over $ k\in[0,T]$, and using~\eqref{zo_pb:step:eta2-sm}  yield
\begin{align}\label{zo_pb:zerosg:thm-sg-sm-equ4p}
&n(\mathbb{E}[f(\bar{x}_{T+1})]-f^*)=\mathbb{E}[W_{4,T+1}]\nonumber\\
&\le W_{4,0}+\frac{2\sqrt{n}}{\sqrt{pT}} L_f^2\sum_{k=0}^{T}\|\bsx_k\|^2_{\bsK}+n\mathcal{O}(p)\eta^2\frac{T+1}{T}\nonumber\\
&~~~+\mathcal{O}(np)\eta\delta^2_k\sqrt{\frac{T+1}{T}}.
\end{align}

Noting that $W_{4,0}=\mathcal{O}(n)$ and $\sqrt{n}n/\sqrt{pT}<1$ due to $T> n^3/p$, from \eqref{zo_pb:zerosg:thm-sg-sm-equ3.1p} and~\eqref{zo_pb:zerosg:thm-sg-sm-equ4p}, we have $\mathbb{E}[f(\bar{x}_{T+1})]-f^*=\mathcal{O}(1)$, which gives~\eqref{zerosg-p:thm-sg-sm-equ4-ZOOM-PB}.

\section{Numerical Examples}

\subsection{Black-box Binary Classification}\label{exp:bbc}

We first consider the same non-linear least square problem in \cite{xu2020second, liu2019signsgd, liu2018zeroth}, i.e., problem with $f_i(\mathbf x) = \left ( y_i - \phi(\mathbf x; \mathbf a_i) \right )^2$ for $i\in[n]$, where $\phi(\mathbf x; \mathbf a_i) = \frac{1}{1 + e^{-\mathbf a_i^{T} \mathbf x}}$. For preparing the synthetic dataset, we randomly draw samples $\mathbf a_i$ from $\mathcal{N}(\mathbf 0, \mathbf{I} )$, and we set the optimal vector $\mathbf{x_{opt}} = \mathbf{1}$, the label is $y_i = 1$ if $\phi(\mathbf x_{opt}; \mathbf a_i) \geq 0.5$ and $0$ otherwise. The training set has $2000$ samples and the test set has $200$ samples. We set the dimension $d$ of $\mathbf a_i$ as $100$, batchsize is $1$, and the total iteration number as $10000$. As suggested in the work \cite{liu2019signsgd}, the smooth parameter $\delta = \frac{10}{\sqrt{Td}}$.

We compare the proposed algorithms with state-of-the-art dencentralized ZO algorithms: ZO-GDA \cite{tang2020distributedzero}, ZONE-M \cite{hajinezhad2019zone}, ZODPDA, ZODPA \cite{yi2021zerothorder}, ZODIAC~\cite{zhang2021convergence}, and ZODIAC-PB~\cite{zhang2021zo}.
The communication topology of $10$ agents is generated randomly following the Erd\H{o}s - R\' enyi model with probability of $0.4$.

\begin{figure}[!ht]
\centering
  \includegraphics[width=0.5\textwidth]{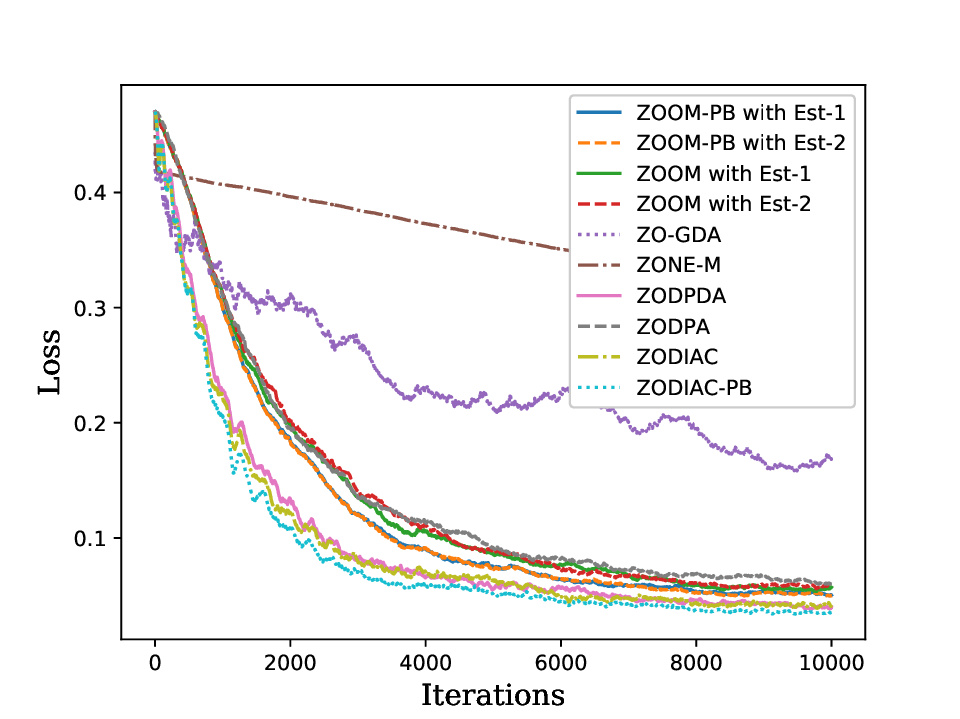}
  \caption{Performance comparison of training loss.}
  \label{nonconvex:fig:bb}
\end{figure}

The comparison of training loss with respect to each algorithm is shown in Figure~\ref{nonconvex:fig:bb}, and the accuracy of each algorithm is given in Table~\ref{tab:acc}, which shows that the proposed algorithms can obtain similar accuracy with other existing algorithms.
From Figure~\ref{nonconvex:fig:bb}, we can see that ZOOM converges similarly as ZODPA which has the convergence rate of $\mathcal{O}(\sqrt{p}/\sqrt{nT})$, and ZOOM-PB~($\gamma = 0.7$) converges faster than ZOOM.
However, we can easily see that ZODPDA, ZODIAC and ZODIAC-PB have better convergence results in this case because all of them communicate both primal and dual variables in each iteration.

\begin{table}[ht!]
\caption{Accuracy}
\label{tab:acc}
\begin{center}
\begin{tabular}{cc}
\multicolumn{1}{c}{Algorithm}  & Accuracy($\%$)
\\ \hline
ZOOM-PB with Est-1   & \bf{97.0}  \\
ZOOM-PB with Est-2   & \bf{96.0}  \\
ZOOM with Est-1  & \bf{96.0}  \\
ZOOM with Est-2  & \bf{96.5} \\
ZO-GDA       & 92.5  \\
ZONE-M       & 89.5  \\
ZODPDA       & 97.5  \\
ZODPA        & 95.5  \\
ZODIAC       & 97.0  \\
ZODIAC-PB    & 97.5  \\
\hline
\end{tabular}
\end{center}
\end{table}

\subsection{Generation of Adversarial Examples from Black-box DNNs}\label{exp:dnn}

In this section, we consider the benchmark example of generation of adversarial examples from black-box DNNs in ZO optimization literature \cite{chen2017zoo, liu2019signsgd, liu2018zeroth}.
In image classification tasks, convolutional neural networks are vulnerable to adversarial examples \cite{goodfellow2014explaining} even under small perturbations, which leads to misclassifications.
Considering the setting of zeroth-order attacks \cite{carlini2017towards}, the model is hidden and no gradient information is available. We treat this task of generating adversarial examples as a zeroth-order optimization problem.
Formally, the loss function\footnote{The loss function violates Assumption~\ref{zo_pb:ass:zeroth-smooth}, however, this experiment is considered as a benchmark in zeroth-order optimization literature. Moreover, this motivates us to investigate nonsmooth problems in the future.} is given as in~\eqref{zo_pb:exp:loss_gan} 
\begin{equation}\label{zo_pb:exp:loss_gan}
\begin{split}
f_i(\mathbf x) = & c \cdot \max \{ F_{y_i} (0.5 \cdot \tanh ( \tanh^{-1} 2 \mathbf  a_i + \mathbf x)) \\
&-\max_{j \neq y_i} F_j(0.5 \cdot \tanh ( \tanh^{-1} 2 \mathbf a_i + \mathbf x)) , 0 \} \\
& + \|0.5 \cdot \tanh ( \tanh^{-1} 2 \mathbf a_i + \mathbf x) - \mathbf a_i \|_2^2
\end{split}
\end{equation}
where $(\mathbf a_i,y_i)$ denotes the pair of the $i$th natural image $\mathbf a_i$ and its original class label $y_i$. The output of function $F(\mathbf z)=[F_1(\mathbf z),\ldots,F_N(\mathbf z)]$ is the well-trained model prediction of the input $\mathbf z$ in all $N$ image classes.
The well-trained DNN model\footnote{\url{https://github.com/carlini/nn_robust_attacks}} on MNIST handwritten has $99.4\%$ test accuracy on natural examples \cite{liu2018zeroth}.
The purpose of this experiment is to generate \textit{false} examples to attack the DNN model in order to have a wrong prediction, i.e. if feeding an original image with label $1$, the DNN predicts it as $1$, however after generating the \textit{false} example based on the original $1$, the DNN should make a wrong prediction.
We compare the proposed algorithms with the same existing algorithms in Section~\ref{exp:bbc}.

\begin{figure}[!ht]
\centering
  \includegraphics[width=0.5\textwidth]{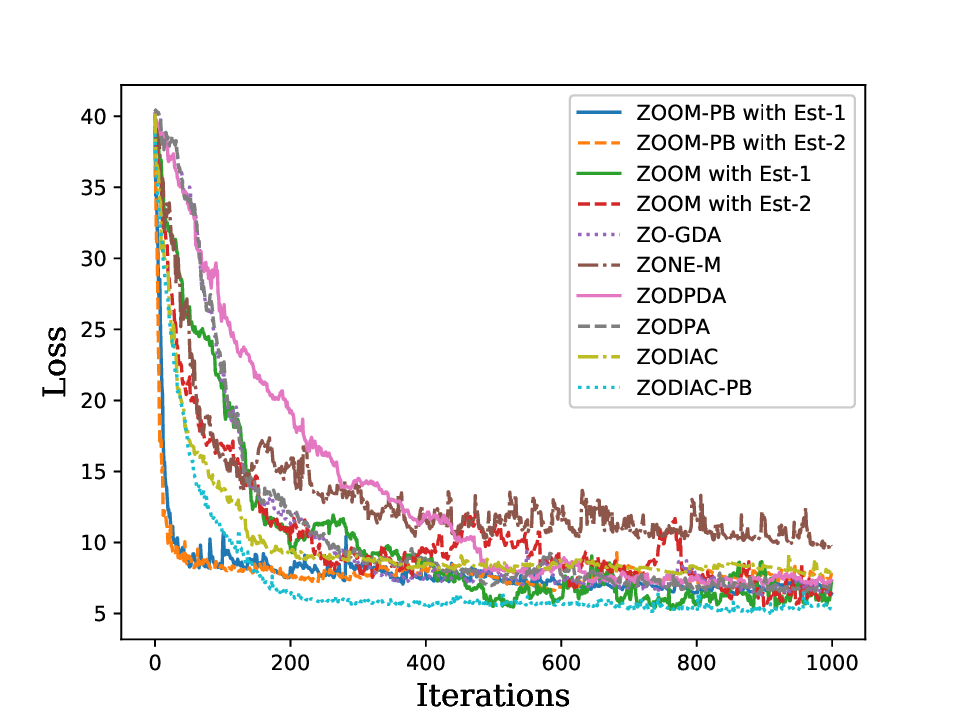}
  \caption{Performance comparison of training loss.}
  \label{nonconvex:fig:dnn:loss}
\end{figure}

\begin{figure}[!ht]
\centering
  \includegraphics[width=0.25\textwidth]{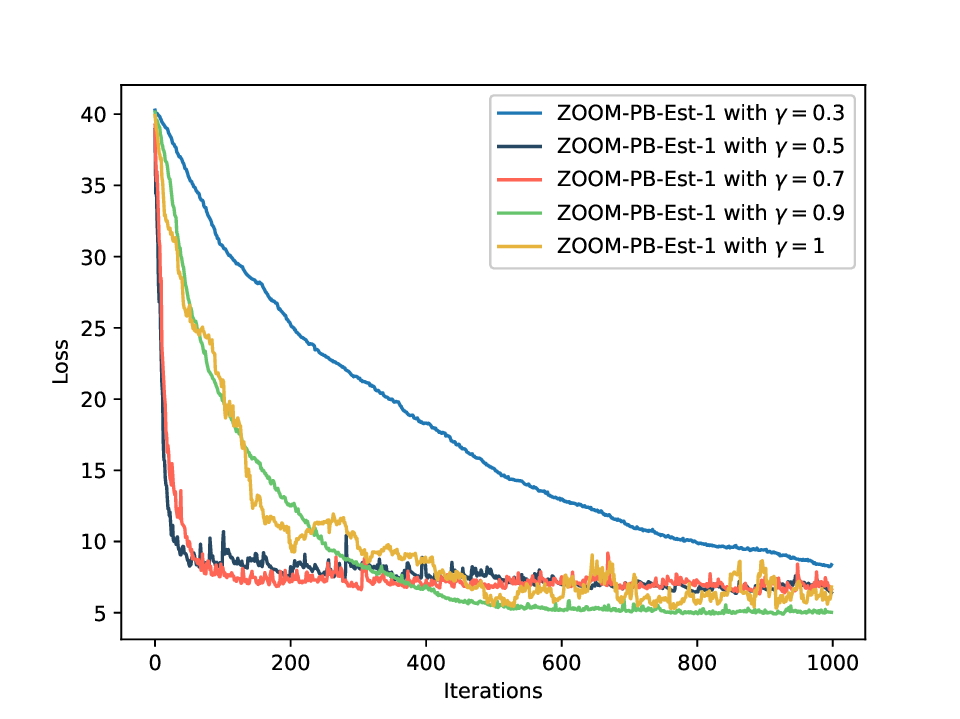}~\includegraphics[width=0.25\textwidth]{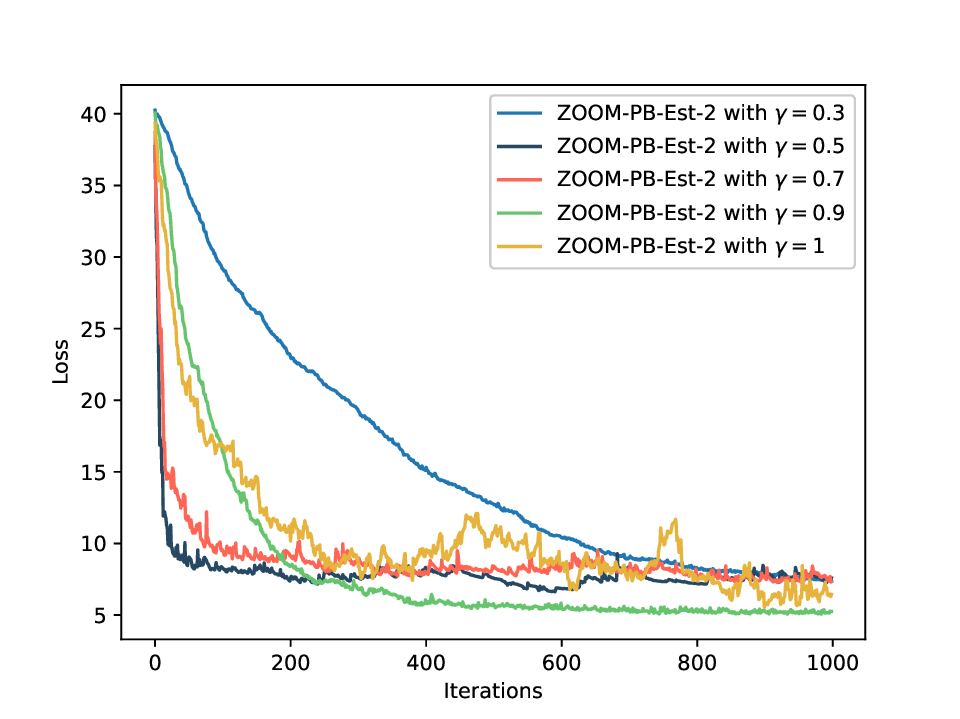}
  \caption{Robustness of $\gamma$ for two gradient estimators.}
  \label{nonconvex:fig:dnn:rb}
  
\end{figure}

The comparison of training loss with respect to each algorithm is shown in Figure~\ref{nonconvex:fig:dnn:loss}, the comparison of distortion is summarized in Table~\ref{tab:dist}, and the generated examples from proposed algorithms are demonstrated in Table~\ref{table:digit4}.
From Figure~\ref{nonconvex:fig:dnn:loss}, we can see that in this non-smooth case, the proposed algorithm ZOOM also has a similar convergence behavior with existing algorithms with convergence rate of $\mathcal{O}(\sqrt{p}/\sqrt{nT})$, and ZOOM-PB~($\gamma = 0.7$) converges fastest among all implemented algorithms.
Moreover, from Table~\ref{tab:dist}, we can conclude that ZOOM and ZOOM-PB have better results in terms of distortion.
Additionally, we test the robustness of $\gamma$ under the same experimental setup.
In Figure~\ref{nonconvex:fig:dnn:rb}, we show the training loss with $\gamma \in \{0.3, 0.5, 0.7, 0.9, 1\}$ for both gradient estimators of ZOOM-PB.
We can conclude that ZOOM-PB is robust to $\gamma \in [0.5, 1]$ and the demonstration matches our theoretical results.

\begin{table}[ht!]
\caption{Distortion}
\label{tab:dist}
\begin{center}
\begin{tabular}{cc}
\multicolumn{1}{c}{Algorithm}  & $l_2$ Distortion
\\ \hline
ZOOM-PB with Est-1   & \bf{4.82}  \\
ZOOM-PB with Est-2   & \bf{4.95}  \\
ZOOM with Est-1  & \bf{5.82}  \\
ZOOM with Est-2  & \bf{6.62} \\
ZO-GDA       & 7.23  \\
ZONE-M       & 9.96  \\
ZODPDA       & 6.44  \\
ZODPA        & 5.77  \\
ZODIAC       & 7.18  \\
ZODIAC-PB    & 5.23  \\
\hline
\end{tabular}
\end{center}
\end{table}
\begin{table}[!ht] \caption{Demonstration of generated adversarial examples from a black-box DNN on MNIST: digit class ``4''.} \label{table:digit4}
  \begin{center}
  \begin{adjustbox}{max width=0.5\textwidth}
  \begin{tabular}
      {cccccc}
      \hline
      	Image ID & 4 & 6 & 19 & 24 & 27 \\
      \hline &&&&& \vspace{-0.2cm} \\
      	Original &
        \parbox[c]{2.2em}{\includegraphics[width=0.4in]{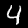}} &
        \parbox[c]{2.2em}{\includegraphics[width=0.4in]{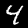}} &
        \parbox[c]{2.2em}{\includegraphics[width=0.4in]{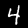}} &
        \parbox[c]{2.2em}{\includegraphics[width=0.4in]{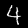}} &
        \parbox[c]{2.2em}{\includegraphics[width=0.4in]{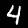}}  \\
        Classified as & 4 & 4 & 4 & 4 & 4  \\
      \hline &&&&& \vspace{-0.2cm} \\
      	ZOOM with Est-1 &
        \parbox[c]{2.2em}{\includegraphics[width=0.4in]{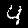}} &
        \parbox[c]{2.2em}{\includegraphics[width=0.4in]{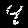}} &
        \parbox[c]{2.2em}{\includegraphics[width=0.4in]{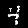}} &
        \parbox[c]{2.2em}{\includegraphics[width=0.4in]{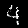}} &
        \parbox[c]{2.2em}{\includegraphics[width=0.4in]{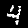}}  \\
        Classified as & 9 & 8 & 2 & 8 & 2  \\
       \hline &&&&& \vspace{-0.2cm} \\
      	ZOOM with Est-2 &
        \parbox[c]{2.2em}{\includegraphics[width=0.4in]{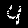}} &
        \parbox[c]{2.2em}{\includegraphics[width=0.4in]{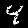}} &
        \parbox[c]{2.2em}{\includegraphics[width=0.4in]{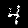}} &
        \parbox[c]{2.2em}{\includegraphics[width=0.4in]{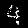}} &
        \parbox[c]{2.2em}{\includegraphics[width=0.4in]{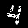}}   \\
        Classified as & 9 & 8 & 2 & 3 & 2  \\
       \hline &&&&& \vspace{-0.2cm} \\
      	ZOOM-PB with Est-1 &
        \parbox[c]{2.2em}{\includegraphics[width=0.4in]{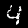}} &
        \parbox[c]{2.2em}{\includegraphics[width=0.4in]{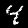}} &
        \parbox[c]{2.2em}{\includegraphics[width=0.4in]{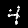}} &
        \parbox[c]{2.2em}{\includegraphics[width=0.4in]{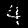}} &
        \parbox[c]{2.2em}{\includegraphics[width=0.4in]{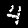}}  \\
        Classified as & 9 & 8 & 2 & 2 & 2  \\
       \hline &&&&& \vspace{-0.2cm} \\
      	ZOOM-PB with Est-2 &
        \parbox[c]{2.2em}{\includegraphics[width=0.4in]{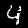}} &
        \parbox[c]{2.2em}{\includegraphics[width=0.4in]{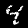}} &
        \parbox[c]{2.2em}{\includegraphics[width=0.4in]{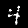}} &
        \parbox[c]{2.2em}{\includegraphics[width=0.4in]{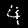}} &
        \parbox[c]{2.2em}{\includegraphics[width=0.4in]{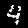}}  \\
        Classified as & 9 & 1 & 7 & 9 & 9  \\
        
       \hline
  \end{tabular}
  \end{adjustbox}
  \end{center}
\end{table}

\section{Conclusion}

In this paper, we investigated the ZO stochastic decentralized non-convex optimization problems and proposed ZOOM and its acceleration variant ZOOM-PB.
We demonstrated that the proposed algorithms can achieve the convergence rate of $\mathcal{O}(\sqrt{p}/\sqrt{nT})$ for general non-convex cost functions.
Additionally, we illustrated the efficacy of considered algorithms through benchmark examples on a large-scale multi-agent topology in comparison with the existing state-of-the-art centralized and distributed ZO algorithms.

\bibliographystyle{IEEEtran}
\bibliography{bibliography}

\appendix\label{zo_pb:sec:app}

\section*{Appendix}
\section{Proof of Lemma~4}\label{zo_pb:proof:lemma4}
\begin{proof}

Consider the following Lyapunov candidate function 
\begin{align}
W_{k}= &\underbrace{\frac{1}{2}\|\bsx_{k}\|^2_{\bsK}}_{W_{1, k}} + \underbrace{n(f(\bar{x}_k)-f^*)}_{W_{2, k}}
\end{align}
where $\bsQ=R\Lambda^{-1}_1R^{\top}\otimes {\bm I}_p$.
Additionally, we denote $g^s_{i,k}=\nabla \mathbb{E}[f_i(x+\delta_{i,k} e_{i})]$, $\bsg^s_k=\col(g^s_{1,k},\dots,g^s_{n,k})$, $\bar{\bsg}^s_k=\bsH\bsg^s_k$, $\bar{g}^e_k=\frac{1}{n}({\bm 1}_n^\top\otimes{\bm I}_p)\bsg^e_k$, and $\bar{\bsg}^e_k={\bm 1}_n\otimes\bar{g}^e_k=\bsH\bsg^e_k$.

(i) We have
\begin{align}
&\mathbb{E}[W_{1,k+1}]
=\mathbb{E}\Big[\frac{1}{2}\|\bm{x}_{k+1} \|^2_{\bsK}\Big]\nonumber\\
&=\mathbb{E}\Big[\frac{1}{2}\|\bm{x}_k-(\alpha\bsL\bm{x}_k+\eta\sigma(\bsg^e_k, \gamma)) \|^2_{\bsK}\Big]\nonumber\\
&\overset{\text{(a)}}{=}\mathbb{E}\Big[\frac{1}{2}\|\bm{x}_k\|^2_{\bsK}-\alpha\|\bsx_k\|^2_{\bsL}
+\frac{1}{2}\alpha^2\|\bsx_k\|^2_{\bsL^2}
\nonumber\\
&~~~-\eta\bsx^\top_k({\bm I}_{np}-\alpha\bsL)\bsK\sigma(\bsg^e_k, \gamma)\nonumber\\
&~~~+\frac{1}{2}\eta^2\Big\|\sigma(\bsg^e_k, \gamma)\Big\|^2_{\bsK}\Big]\nonumber\\
&\overset{\text{(b)}}{=}W_{1,k}-\|\bsx_k\|^2_{\alpha\bsL
-\frac{1}{2}\alpha^2\bsL^2}\nonumber\\
&~~~-\eta\bsx^\top_k({\bm I}_{np}-\alpha\bsL)\bsK\bsg^s_k\nonumber\\
&~~~+\frac{1}{2}\eta^2\mathbb{E}\Big[\Big\|\bsg_k^0
+\sigma(\bsg^e_k, \gamma)-\bsg_k^0\Big\|^2_{\bsK}\Big]\nonumber\\
&\overset{\text{(c)}}{\le} W_{1,k}-\|\bsx_k\|^2_{\alpha\bsL
-\frac{1}{2}\alpha^2\bsL^2}
-\eta\bsx^\top_k\bsK\bsg_k^0\nonumber\\
&~~~+\frac{1}{2}\|\bm{x}_k\|^2_{\bsK}
+\frac{1}{2}\|\bsg^s_k-\bsg_k^0\|^2\nonumber\\
&~~~+\alpha^2\|\bm{x}_k\|^2_{\bsL^2}
+\frac{1}{2}\beta^2\Big\|\bsg_k^0\Big\|^2_{\bsK}
+\frac{1}{2}\|\bsg^s_k-\bsg_k^0\|^2\nonumber\\
&~~~+\eta^2\Big\|\bsg_k^0\Big\|^2_{\bsK}
+\mathbb{E}[\|\sigma(\bsg^e_k, \gamma)-\bsg_k^0\|^2]\nonumber\\
&\overset{\text{(d)}}{\le} W_{1,k}-\|\bsx_k\|^2_{\alpha\bsL-\frac{1}{2}\bsK
-\alpha^2\bsL^2-\eta(1+5\eta)L_f^2\bsK}
\nonumber\\
&~~~-\beta\bsx^\top_k\bsK\Big(\bsg_k^0\Big)
+\Big\|\bsg_k^0\Big\|^2_{\frac{3}{2}\beta^2\bsK}\nonumber\\
&~~~+nL_f^2\Big[\frac{p}{2}+4\Big]\delta_k^2+2\mathbb{E}[\|\sigma(\bsg^e_k, \gamma)\|^2],\label{zo_pb:zerosg:v1k}
\end{align}
where (a) holds due to Lemma~1 and~2 in~\cite{Yi2018distributed}; (b) holds due to $\mathbb{E}[\bsg^e_k] = \bsg^s_k$ and that $x_{i,k}$ and $v_{i,k}$ are independent of $u_{i,k}$ and $\xi_{i,k}$; (c) holds due to the Cauchy--Schwarz inequality and $\rho(\bsK)=1$; (d) holds due to $\|\bsg^s_k-\bsg_k^0\|^2\le 2L_f^2\|\bsx_{k}\|^2_{\bsK}+\frac{np}{2}L_f^2\delta_k^2$ and $\mathbb{E}[\|\bsg_k^0-\sigma(\bsg^e_k, \gamma)\|^2]\le 4L_f^2\|\bsx_{k}\|^2_{\bsK}+4nL_f^2\delta_k^2+2\mathbb{E}[\|\sigma(\bsg^e_k, \gamma)\|^2]$.

(ii) We have
\begin{align}
&\mathbb{E}[W_{2,k+1}]=\mathbb{E}[n(f(\bar{x}_{k+1})-f^*)]
=\mathbb{E}[\tilde{f}(\bar{\bsx}_{k+1})-nf^*]\nonumber\\
&=\mathbb{E}[\tilde{f}(\bar{\bsx}_k)-nf^*+\tilde{f}(\bar{\bsx}_{k+1})
-\tilde{f}(\bar{\bsx}_k)]\nonumber\\
&\overset{\text{(e)}}{\le}\mathbb{E}[\tilde{f}(\bar{\bsx}_k)-nf^*
-\eta(\bar{\bsg}_{k}^e)^\top\bsg^0_k
+\frac{1}{2}\eta^2L_f\|\bar{\bsg}_{k}^e\|^2]\nonumber\\
&\overset{\text{(f)}}{=}W_{2,k}
-\eta(\bar{\bsg}_{k}^s)^\top\bar{\bsg}^0_k
+\frac{1}{2}\eta^2L_f\mathbb{E}[\|\bar{\bsg}_{k}^e\|^2]\nonumber\\
&\overset{\text{(g)}}{=}W_{2,k}
-\frac{1}{2}\eta(\bar{\bsg}_{k}^s)^\top(\bar{\bsg}^s_k+\bar{\bsg}^0_k-\bar{\bsg}^s_k)\nonumber\\
&~~~-\frac{1}{2}\eta(\bar{\bsg}^s_{k}-\bar{\bsg}^0_k+\bar{\bsg}^0_k)^\top\bar{\bsg}^0_k
+\frac{1}{2}\eta^2L_f\mathbb{E}[\|\bar{\bsg}_{k}^e\|^2]\nonumber\\
&\overset{\text{(h)}}{\le} W_{2,k}-\frac{1}{4}\eta(\|\bar{\bsg}^s_{k}\|^2
-\|\bar{\bsg}^0_k-\bar{\bsg}^s_k\|^2+\|\bar{\bsg}_{k}^0\|^2\nonumber\\
&~~~-\|\bar{\bsg}^0_k-\bar{\bsg}^s_k\|^2)
+\frac{1}{2}\eta^2L_f\mathbb{E}[\|\bar{\bsg}_{k}^e\|^2]\nonumber\\
&\overset{\text{(i)}}{\le} W_{2,k}-\frac{1}{4}\eta\|\bar{\bsg}^s_{k}\|^2
+\|\bsx_k\|^2_{\eta L_f^2\bsK}\nonumber\\
&~~~+\frac{np}{4}\eta L_f^2\delta^2_k-\frac{1}{4}\eta\|\bar{\bsg}_{k}^0\|^2
+\frac{1}{2}\eta^2L_f\mathbb{E}[\|\bar{\bsg}^e_{k}\|^2],\label{zerosg:v4k}
\end{align}
where (e) holds since that $\tilde{f}$ is smooth; (f) holds due to $\mathbb{E}[\bsg^e_k] = \bsg^s_k$, $x_{i,k}$ and $v_{i,k}$ are independent; (g) holds due to $(\bar{\bsg}_{k}^s)^\top\bsg^0_k=(\bsg_{k}^s)^\top\bsH\bsg^0_k=(\bsg_{k}^s)^\top\bsH\bsH\bsg^0_k
=(\bar{\bsg}_{k}^s)^\top\bar{\bsg}^0_k$; (h) holds due to the Cauchy--Schwarz inequality; and (i) holds due to  $\|\bsg^s_k-\bsg_k^0\|^2\le 2L_f^2\|\bsx_{k}\|^2_{\bsK}+\frac{np}{2}L_f^2\delta_k^2$.

(iii) 
Define $W_{k+1} = W_{1, k+1} + W_{2, k+1}$ and then we have the following inequality holds.

\begin{align}
&\mathbb{E}[W_{k+1}]\nonumber\\
&\le W_{k}-\|\bsx_k\|^2_{\alpha\bsL-\frac{1}{2}\bsK
-\alpha^2\bsL^2-\eta(1+5\eta)L_f^2\bsK}
\nonumber\\
&~~~-\beta\bsx^\top_k\bsK\Big(\bsg_k^0\Big)
+\Big\|\bsg_k^0\Big\|^2_{\frac{3}{2}\beta^2\bsK}\nonumber\\
&~~~+nL_f^2\Big[\frac{p}{2}+4\Big]\delta_k^2+2\mathbb{E}[\|\sigma(\bsg^e_k, \gamma)\|^2]\nonumber\\
&~~~-\frac{1}{4}\eta\|\bar{\bsg}^s_{k}\|^2
+\|\bsx_k\|^2_{\eta L_f^2\bsK}\nonumber\\
&~~~+\frac{np}{4}\eta L_f^2\delta^2_k-\frac{1}{4}\eta\|\bar{\bsg}_{k}^0\|^2
+\frac{1}{2}\eta^2L_f\mathbb{E}[\|\bar{\bsg}^e_{k}\|^2]\\
&\overset{\text{(j)}}{\le} W_{k}-\|\bsx_k\|^2_{\bsM_{1}-\bsM_{2}-b_{1}\bsK}-\Big\|\bsg_{k}^0\Big\|^2_{b_{2}\bsK}\nonumber\\
&~~~-\eta\Big(\frac{1}{4}-6c_1(p-1)\Big)\|\bar{\bsg}^0_{k}\|^2\nonumber\\
&~~~+\underbrace{c_1\Big[ 6(p-1)\sigma^2_2+(\frac{3}{n_c}+2)p\sigma^2_1 \Big]}_{\mathcal{O}(np)\eta^2}+\underbrace{c_{3}\eta\delta_k^2L_f^2}_{\mathcal{O}(np^2)\eta\delta_k^2},
\label{zo_pb:zerosg:vkLya}
\end{align}
where (j) holds due to \eqref{zo_pb:zerosg:rand-grad-esti2}, \eqref{zo_pb:zerosg:rand-grad-esti4}, and
\begin{align*}
\bsM_{1}&=\alpha\bsL-\Big(1+3L_f^2+6L_f^4\kappa_1(p-1)\Big)\bsK,\\
\bsM_{2}&=\bsL+2\alpha^2\bsL^2+8L_f^2\bsK\\
&\quad+6(p-1)\Big(3+\frac{1}{2}L_f+{2L_f^{2}}\kappa_1 + \frac{L_f^2}{2}\Big)\bsK,\\
\kappa_1&>\frac{1}{\rho_2(L)}+1, \\
\kappa_2&\in\Big(0,\min\{\frac{(\kappa_1-1)\rho_2(L)-1}{\rho(L)+(2\kappa_1^2+1)\rho(L^2)+1}, \frac{1}{5}\}\Big)\\
\kappa_3&=\frac{1}{\rho_2(L)}+\kappa_1+1,\\
b_{1}
&=6p\kappa_3L_f^4{\eta}+12p(\kappa_3+1)L_f^4{\eta^2},\\
b_{2}&=\frac{1}{2}\eta(2\beta-\kappa_3)-2.5\kappa_2^2,\\
c_1&= \Big(3+\frac{1}{2}L_f+{2L_f^{2}}\kappa_1 + {L_f^2}{2}\Big)n\eta^2 + {L_f^2\kappa_1}n\eta,\\
c_2 &=\frac{3}{4}pn+\eta n(\frac{p}{2} +6) +\frac{p^2L_f^2}{2}\kappa_1,\\
c_3 &=c_2 + \Big(c_1 - {L_f^2\kappa_1}n\eta\Big)p^2\eta,\\
\epsilon&=\frac{1}{2}\alpha\rho_2(L)-\alpha^2\rho(L^2).
\end{align*}

Consider $p\geq 1$, $\alpha \in (0, \frac{\rho_2(L)}{2\rho(L^2)})$, and \\
$\eta \in (0,\min \{\frac{2\alpha\rho_2(L)-4\alpha^2\rho(L^2)}{9L_f^{2}}, \frac{\alpha\rho_2(L)}{48p[(1+2\alpha\rho_2(L))+\alpha\rho_2(L)L_f]}\}\big]$, we have~\eqref{zo_pb:zerosg:sgproof-vkLya2T}, and~\eqref{zo_pb:zerosg:v4kspeed}.
\end{proof}

\end{document}